\documentclass[11pt,a4paper]{article}

\usepackage{latexsym}
\usepackage{amsmath}
\usepackage{epsfig}
\usepackage{times}
\usepackage{amssymb}
\newtheorem{thr}{\quad Theorem}
\newtheorem{lem}{\quad Lemma}

\newtheorem{exa}{\quad Example}
\newcommand{\Z}{\mathbb{Z}}

\title{A method to construct generalized balanced tournament designs}

\author{Songchol Kim and Changil Rim \\ \\
\small\textit{Department of Mathematics, } \\
\small\textit{Kim Il Sung University,  Pyongyang D.P.R.Korea}\\
\small\textit {email: songchol$\_$kim@yahoo.com}}

\begin{document}
\maketitle
%\hline

\begin{abstract}
A generalized balanced tournament design, or a GBTD$(k, m)$ in short, is a $(km, k, k-1)-$BIBD defined on a $km$-set $V$ . Its blocks can be arranged into an $m\times (km-1)$ array in such a way that (1) every element of $V$ is contained in exactly one cell of each column, and (2) every element of $V$ is contained in at most $k$ cells of each row. In this paper, we present a new construction for GBTDs and show that a GBTD$(p,p)$ exists for any prime number $p\geq 3$.

\bigskip



\end{abstract}

%\hline
%
%-----------------------1. Introduction-----------------
%
\section{Introduction}

Let $X$ be a set of $v$ elements (called points), and $A$ a collection of subsets (called blocks) of $X$. The ordered pair $(X,A)$ is known as a $(v, k, \lambda)$ {\it balanced incomplete block design}(BIBD), or a $(v, k, \lambda )-$BIBD, if every pair of distinct points occurs in precisely $\lambda$  blocks.\\

The number of blocks of a $(v, k, \lambda )-$BIBD is $\lambda v(v-1)/k(k-1)$(\cite{colbourn}). Hence, $(km, k, k-1)-$BIBD has $m(km-1)$ blocks.\\

Following Lamken \cite{lamken2}, a $(km, k, k-1)-$BIBD $(X,A)$ is called a {\it generalized balanced tournament design} (GBTD), or a GBTD$(k, m)$ in short, when the $m(km-1)$ blocks of $A$ can be arranged into an $m\times (km-1)$ array in such a way that

(1) every point of $X$ occurs in exactly once in each column;

(2) every point of $X$ is contained in at most $k$ cells of each row.\\

By definition, any GBTD can be identified with its corresponding array of blocks. In what follows, we will make no difference between a GBTD and its corresponding array of blocks. The following lemma is straightforward, which was stated in \cite{lamken2}.\\

%
%-----------------Lemma 1--------------------------
%
\begin{lem}
Every point of a GBTD$(k, m)$ occurs $k$ times in $(m-1)$ rows and $(k-1)$ times in the remaining row.
\end{lem}

Now consider a GBTD$(k, m)$, $R$, over $X$. A point contained in only $k-1$ cells of row $i$ of $R$ is called a deficient point of row $i$. It is easily seen that each row of $R$ contains exactly $k$ deficient points which are pairwise distinct. The $k-$tuple consisting of the $k$ deficient points is referred to as the {\it deficient} $k-${\it tuple} of row $i$. As an immediate consequence of Lemma 1, we have the following.

%
%-----------------Lemma 2--------------------------
%
\begin{lem}
The deficient k-tuples of a GBTD(k, m) partition its point set(see [2]).
\end{lem}

It is known (see \cite{colbourn,lamken1}) that for all positive $m{\not =} 2$, both a GBTD$(2,m)$ and a GBTD$(3,m)$ exist, whilst there does not exist a GBTD$(k, 2)$ for all positive $k\geq 2$.\\
It is also known (see \cite{jianxing}) that for any integer $m\geq 5$, a GBTD$(4,m)$ exists with at most eight possible exceptions of $m\in $\{28, 32, 33, 34, 37, 38, 39, 44\}.\\
To our knowledge, very little is known about the existence of a GBTD$(k, m)$ with $k\geq 5$. In this paper we consider the existence of a GBTD$(p,p)$. We show in Sect. 3 that a GBTD$(p,p)$ exists for any prime number $p\geq 3$.

Section 3 contains a new construction for GBTDs. Our construction mainly use the matrix over the field $\Z_p$. It is presented in Sect. 2. A link between a GBTD$(k, m)$ and a near constant composition code is mentioned in \cite{jianxing}. The derived code from GBTD$(k, m)$��s is optimal in the sense of its size.

%\hline
%
%-----------------------Section 2-----------------
%
\section{GBTD$(p,p)$'s and square matrices of order $p^2$ over the field $\Z_p$($p\geq 3$ is a prime number)}

Let $R=(r_{ij})$ be a GBTD$(p,p)$(we index the rows by the elements of $\Z_p$, and the columns by the elements of $\{1,2,\cdots,p^2-1\}$). Now we construct a $p^2$ by $p^2$ matrix $M'=(m_{ij})$ over the field $\Z_p$ from $R$.\\

\indent Let $\{1,2,\cdots,p^2\}$ be a point set of $R$.

$R$ contains every point exactly once in each column.

For all $i$ and $j$ with $1\leq i\leq p^2-1$,$1\leq j\leq p^2$, we set $m_{i+1,j}=k$ if $j$ is contained $r_{ki}$.

From the lemma 1, for all $j$ with $1\leq j\leq p^2$, the multiset $\{m_{2j},m_{3j},m_{4j},\cdots,m_{pj}\}$ contains $p-1$ elements of $\Z_p$ $p$ times respectively, and remaining element $p-1$ times(we denote this element by $d_j$).

We set $m_{1j}$ to $d_j$.

Clearly, the set $\{j:m_{1j}=i\}$ is a deficient $k-$tuple of row $i$ for any $i$ of $\Z_p$.\\

%
%-----------------------Example 1-----------------
%
\begin{exa}
Let
\begin{center}
\begin{tabular}{c|c|c|c|c|c|c|c|c|c|}
\cline{2-10}
 & \textup {129} & \textup {349} & \textup {569} & \textup {145} & \textup {357} & \textup {178} & \textup {238} & \textup {267} & {\bf 468} \\
\cline{2-10}
R= & \textup {357} & \textup {167} & \textup {138} & \textup {236} & \textup {468} & \textup {245} & \textup {749} & \textup {589} & {\bf 129} \\
\cline{2-10}
 & \textup {468} & \textup {258} & \textup {247} & \textup {789} & \textup {129} & \textup {369} & \textup {165} & \textup {134} & {\bf 357} \\
\cline{2-10}
\end{tabular}
\end{center}

In the previous figure, the bold elements present deficient tuples. It can be easily seen that $R$ is a GBTD(3,3). From $R$, the following matrix is obtained.

\begin{center}
\begin{tabular}{c|c|c|c|c|c|c|c|c|c|}
\cline{2-10}
 & {\textup 1} & {\textup 1} & {\textup 2} & {\textup 0} & {\textup 2} & {\textup 0} & {\textup 2} & {\textup 0} & {\textup 1} \\
\cline{2-10}
 & {\textup 0} & {\textup 0} & {\textup 1} & {\textup 2} & {\textup 1} & {\textup 2} & {\textup 1} & {\textup 2} & {\textup 0} \\
\cline{2-10}
 & {\textup 1} & {\textup 2} & {\textup 0} & {\textup 0} & {\textup 2} & {\textup 1} & {\textup 1} & {\textup 2} & {\textup 0} \\
\cline{2-10}
 & {\textup 1} & {\textup 2} & {\textup 1} & {\textup 2} & {\textup 0} & {\textup 0} & {\textup 2} & {\textup 1} & {\textup 0} \\
\cline{2-10}
$M'$= & {\textup 0} & {\textup 1} & {\textup 1} & {\textup 0} & {\textup 0} & {\textup 1} & {\textup 2} & {\textup 2} & {\textup 2} \\
\cline{2-10}
 & {\textup 2} & {\textup 2} & {\textup 0} & {\textup 1} & {\textup 0} & {\textup 1} & {\textup 0} & {\textup 1} & {\textup 2} \\
\cline{2-10}
 & {\textup 0} & {\textup 1} & {\textup 2} & {\textup 1} & {\textup 1} & {\textup 2} & {\textup 0} & {\textup 0} & {\textup 2} \\
\cline{2-10}
 & {\textup 2} & {\textup 0} & {\textup 0} & {\textup 1} & {\textup 2} & {\textup 2} & {\textup 1} & {\textup 0} & {\textup 1} \\
\cline{2-10}
 & {\textup 2} & {\textup 0} & {\textup 2} & {\textup 2} & {\textup 1} & {\textup 0} & {\textup 0} & {\textup 1} & {\textup 1} \\
\cline{2-10}
\end{tabular}
\end{center}
\end{exa}

For any given GBTD$(p,p)$ $R$, the following operations\\
\indent \indent a) interchange of two rows,\\
\indent \indent b) interchange of two columns,\\
\indent \indent c) permutaion on the point set\\
result the equivalent GBTDs.\\

These operations correspond to the following operations on the corresponding matrix over the field $\Z_p$.\\
\indent \indent a) permutation on the $\Z_p$,\\
\indent \indent b) interchange of two rows except the first row,\\
\indent \indent c) interchange of two columns.\\

For a given matrix $M'$, we can permute columns in such a way that the first row becomes
\begin{center}
$(0,0,\cdots,0,1,1,\cdots,1,2,2,\cdots,2,\cdots,p-1,p-1,\cdots,p-1)$.
\end{center}
(where every element of $\Z_p$ is contained exactly $p$ times, respectively)\\

We denote the resultant matrix by $M$, the first $p$ columns of $M$ by $V_0$, next $p$ columns by $V_1,\cdots,$ last $p$ columns by $V_{p-1}$, .i.e.
\begin{center}
$M=[V_0 V_1 \cdots V_{p-1}]$.
\end{center}

From the construction of $M$, it can be seen that $M$ has the following properties;
\indent \indent \textcircled{1} every element of $\Z_p$ is contained exactly $p$ times in each row and each column, respectively,\\
\indent \indent \textcircled{2} for all $i$ with $0\leq i\leq p-1$, any two different columns of $V_i$ contain same entries in exactly $p$ rows,\\
\indent \indent \textcircled{3} for all i and j with $0\leq i<j\leq p-1$, any column of $V_i$ and any one of $V_j$ contain same entries in exactly $p-1$ rows.\\
% \begin{eqnarray}
% && \textcircled{1} \textnormal{ every element of $\Z_p$ is contained exactly $p$ times in each row and each column, respectively,} \\
% && \textcircled{2} \; \textnormal{for all $i$ with $0\leq i\leq p-1$, any two different columns of $V_i$ contain same entries in exactly $p$ rows,}\\
% && \textcircled{3} \; \textnormal{for all i and j with $0\leq i<j\leq p-1$, any column of $V_i$ and any one of $V_j$ contain same entries in exactly $p-1$ rows.}
% \end{eqnarray}

\begin{exa}
In the example\textnormal{1}, any different two columns of $M'$ contain \textnormal{2} same entries in the same rows except the first row, thus for any $j$ and $k$ with $j\not =k$
\begin{center}
$\mid \{ i\geq 2:m_{ij}=m_{ik}\}\mid =2$
\end{center}
For example, the first column and the second column contains same entries respectively in the row 2 and row 6, and the firth column and the fifth column contains same entries respectively in the row 5 and row 7.\\
Every row and every column contains 0,1,2 exactly three times, respectively.
\end{exa}

On the otherhand, it is easily seen that a GBTD$(p,p)$ can be obtained from the matrix which satisfies the previous three conditions.\\

%\hline
%
%-----------------------Section 3-----------------
%
\section{Construction of GBTD$(p,p)$'s}

In this section, we construct $M_p$'s over the field $\Z_p$ which satisfies the three properties mentioned in the section 2.\\

First, we introduce some notations.

We denote the first $p$ rows of $M_p$ by $H^{*}$, and next $p-1$ rows by $H_0$, next $p-1$ rows by $H_1$, next $p-1$ rows by $H_2$,$\cdots$, last $p-1$ rows by $H_{p-1}$.

We denote the first $p$ columns of $M_p$ by $V_0$, and next $p$ columns by $V_1$,$\cdots$, last $p$ columns by $V_{p-1}$.
Thus,
\[ M_p = \left( \begin{array}{c}
H^{*}\\
H_0\\
H_1\\
\cdots\\
H_{p-1}
\end{array} \right) =(V_0,V_1,\cdots,V_{p-1}).
\]

We index the rows of $H^{*}$ by the elements of $\Z_p$.

We index the rows of $H_0,H_1,H_2,\cdots,H_{p-1}$ by the elements of $\Z_p\setminus \{p-1\}$, respectively, and the columns of $V_0,V_1,H_2,\cdots,V_{p-1}$ by the elements of $\Z_p$.

$H_i$ and $V_j$ determine the $(p-1)$ by $p$ matrix. We denote it by $(H_i,V_j)$.

For a given matrix $A$, we denote the row $i$ of $A$ by $(A)_i$, and the column $j$ by $(A)^j$, and the entry in the row $i$ and column $j$ by $(A)^j_i$.

In what follows, operator + and $\times$ are modular operations over the field $\Z_p$.

We denote the vector of order $p$ $(i,i,\cdots,i)$ by $\bar {i}$, and the vector $(i,i+1,i+2,\cdots,i+p-1)$ by $\vec {i}$.\\

We construct a matrix $M_p$ as follows;\\
\begin{center}
$(H^*)_i=(\bar {i},\overline {i+1},\overline {i+2},\cdots,\overline {i+p-1}),\qquad i\in \Z_p$\\
$(H_i,V_j)_l=\overrightarrow {i\times (j+1)+j\times l},\qquad i\in \Z_p,\ \ j,l\in \Z_p\setminus \{p-1\}$\\
\end{center}
\[ (H_i,V_{p-1})_l=\left\{
\begin{array}{l l}
\overrightarrow {-l} & \quad \mbox{if $l\geq i$}\\
\overrightarrow {-l+1} & \quad \mbox{if $l<i$}\\
\end{array} \right. \]

\begin{exa}\ \
\begin{center}
\begin{tabular}{c|c|c|c|c|c|c}
 & $V_0$ & $V_1$ & $V_2$ & $V_3$ & $V_4$ & \\
\cline{2-7}
 & $\bar 0$ & $\bar 1$ & $\bar 2$ & $\bar 3$ & $\bar 4$ & \\
 & $\bar 1$ & $\bar 2$ & $\bar 3$ & $\bar 4$ & $\bar 0$ & \\
 & $\bar 2$ & $\bar 3$ & $\bar 4$ & $\bar 0$ & $\bar 1$ & $H^*$\\
 & $\bar 3$ & $\bar 4$ & $\bar 0$ & $\bar 1$ & $\bar 2$ & \\
 & $\bar 4$ & $\bar 0$ & $\bar 1$ & $\bar 2$ & $\bar 3$ & \\
\cline{2-7}
 & $\vec 0$ & $\vec 0$ & $\vec 0$ & $\vec 0$ & $\vec 0$ & \\
 & $\vec 0$ & $\vec 1$ & $\vec 2$ & $\vec 3$ & $\vec 4$ & \\
 & $\vec 0$ & $\vec 2$ & $\vec 4$ & $\vec 1$ & $\vec 3$ & $H_0$\\
 & $\vec 0$ & $\vec 3$ & $\vec 1$ & $\vec 4$ & $\vec 2$ & \\
\cline{2-7}
 & $\vec 1$ & $\vec 2$ & $\vec 3$ & $\vec 4$ & $\vec 1$ & \\
 & $\vec 1$ & $\vec 3$ & $\vec 0$ & $\vec 2$ & $\vec 4$ & \\
 & $\vec 1$ & $\vec 4$ & $\vec 2$ & $\vec 0$ & $\vec 3$ & $H_1$\\
$M_5=$ & $\vec 1$ & $\vec 0$ & $\vec 4$ & $\vec 3$ & $\vec 2$ & \\
\cline{2-7}
 & $\vec 2$ & $\vec 4$ & $\vec 1$ & $\vec 3$ & $\vec 1$ & \\
 & $\vec 2$ & $\vec 0$ & $\vec 3$ & $\vec 1$ & $\vec 0$ & \\
 & $\vec 2$ & $\vec 1$ & $\vec 0$ & $\vec 4$ & $\vec 3$ & $H_2$\\
 & $\vec 2$ & $\vec 2$ & $\vec 2$ & $\vec 2$ & $\vec 2$ & \\
\cline{2-7}
 & $\vec 3$ & $\vec 1$ & $\vec 4$ & $\vec 2$ & $\vec 1$ & \\
 & $\vec 3$ & $\vec 2$ & $\vec 1$ & $\vec 0$ & $\vec 0$ & \\
 & $\vec 3$ & $\vec 3$ & $\vec 3$ & $\vec 3$ & $\vec 4$ & $H_3$\\
 & $\vec 3$ & $\vec 4$ & $\vec 0$ & $\vec 1$ & $\vec 2$ & \\
\cline{2-7}
 & $\vec 4$ & $\vec 3$ & $\vec 2$ & $\vec 1$ & $\vec 1$ & \\
 & $\vec 4$ & $\vec 4$ & $\vec 4$ & $\vec 4$ & $\vec 0$ & \\
 & $\vec 4$ & $\vec 0$ & $\vec 1$ & $\vec 2$ & $\vec 4$ & $H_4$\\
 & $\vec 4$ & $\vec 1$ & $\vec 3$ & $\vec 0$ & $\vec 3$ & \\
\cline{2-7}

\end{tabular}
\end{center}
\end{exa}

%
%-----------------Lemma 3--------------------------
%
\begin{lem}
For all $m\in \Z_p$ ,the following equation
\begin{equation}
\left\{
\begin{array}{l}
x+y=m\\
0\leq y<p-1
\end{array} \right.
\end{equation}
has $p-1$ solution $(x,y)$'s, and exactly $\frac {p-1}{2}$ solutions of them satisfy $y\geq x$.
\end{lem}
\noindent{\it Proof.}\ \ \ There exists a pair which satisfies (1), for all y with $0\leq y<p-1$. It follows that (1) has $p-1$ solutions.

Now, we prove that $\frac {p-1}{2}$ solutions of them satisfy $y\geq x$.

\noindent (a) If $m=p-1$ then the solutions of (1) are
\[ \{ (1,p-2),(2,p-3),(3,p-4),\cdots,(\frac {p-1}{2},\frac {p-1}{2}),\cdots,(p-1,0)\}\] 
,where the first $\frac {p-1}{2}$ solutions satisfy $y\geq x$.

\noindent (b) If $m<p-1$ then the set of the solutions of the following equations over the field $\Z_p$
\begin{equation}
\left\{
\begin{array}{l}
x+y=m\\
0\leq y<p-1\\
y\geq x
\end{array} \right.
\end{equation}
is sum of two solution sets of the following equations over the nonnegative integers;
\begin{equation}
\left\{
\begin{array}{l}
x+y=m\\
0\leq y\leq m\\
y\geq x
\end{array} \right.
\end{equation}

\begin{equation}
\left\{
\begin{array}{l}
x+y=p+m\\
m<y<p-1\\
y\geq x
\end{array} \right.
\end{equation}

\noindent - In the case that $m$ is even

The solution set of (3) is
\[ \{ (0,m),(1,m-1),\cdots,(\frac {m}{2},\frac {m}{2}),\cdots,(m,0)\}\]
,where only first $\frac {m}{2}+1$ solutions satisfy (2).

The solution set of (4) is
\[ \{ (m+2,p-2),(m+3,p-3),\cdots,(\frac {p+m-1}{2},\frac {p+m+1}{2}),\cdots,(p-1,m+1)\}\]
,where only first $\frac {p+m-1}{2}-(m+1)$ solutions satisfy (2).

Hence, the total number of solutions of (2) is $\frac {p-1}{2}$.

\noindent - In the case that $m$ is odd

The solution set of (3) is
\[ \{ (0,m),(1,m-1),\cdots,(\frac {m-1}{2},\frac {m+1}{2}),\cdots,(m,0)\}\]
,where only first $\frac {m-1}{2}+1$ solutions satisfy (2).

The solution set of (4) is
\[ \{ (m+2,p-2),(m+3,p-3),\cdots,(\frac {p+m}{2},\frac {p+m}{2}),\cdots,(p-1,m+1)\}\]
,where only first $\frac {p+m}{2}-(m+1)$ solutions satisfy (2).

Hence, the total number of solutions of (2) is $\frac {p-1}{2}$.$\Box$\\

\begin{thr}
$M_p$ constructed above satisfies the properties \textcircled{1},\textcircled{2},\textcircled{3} mentioned in sec.2. In other words, a GBTD$(p,p)$ can be obtained from $M_p$.
\end{thr}

\noindent {\it Proof}\ \ \ We can easily see that $M_p$ satifies \textcircled{1} from it's construction.

For all $i\in \Z_p$, any two different columns contain same entries only in the first $p$ rows.

Now we prove that $M_p$ satisfies the property \textcircled{3}.

Let's consider two columns $(V_{j_1})^m$ and $(V_{j_2})^n$, where $0\leq j_1<j_2\leq p-1$,\ $m,n\in \Z_p$.
These two columns don't contain any same entry in the first $p$ rows.
We wish to solve the following equation with respect to $i,l$;\\
\begin{equation}
(H_i,V_{j_1})^m_l=(H_i,V_{j_2})^n_l,\ \ i\in \Z_p,\ ,l\in \Z_p\setminus \{p-1\}.
\end{equation}

\noindent (a) In the case $j_2<p-1$;
\begin{center}
$i\times (j_1+1)+j_1\times l+m=i\times (j_2+1)+j_2\times l+n$\\
$i\times (j_1-j_2)+l\times (j_1-j_2)=n-m$\\
$(i+l)=(n-m)\times (j_1-j_2)^{-1}$\\
\end{center}

Form the lemma 3, it follows that this equation has $(p-1)$ solutions.

Thus, a column of $(V_{j_1})$ and one of $(V_{j_2})$ contain same entries in exactly $(p-1)$ rows.\\

\noindent (b) In the case $j_2=p-1$;

From the definition of $(H_i,V_{p-1})_l$, the solution set of (5) is sum of solution sets of the following two equations;
\begin{equation}
\left\{
\begin{array}{l}
i\times (j_1+1)+j_1\times l+m=-l+n\\
l\geq i
\end{array} \right.
\end{equation}

\begin{equation}
\left\{
\begin{array}{l}
i\times (j_1+1)+j_1\times l+m=-l+1+n\\
l<i
\end{array} \right.
\end{equation}
The equation (6) can be changed into
\begin{equation}
\left\{
\begin{array}{l}
(i+l)\times (j_1+1)=n-m\\
l\geq i
\end{array} \right.
\end{equation}
The equation (7) can be changed into
\begin{equation}
\left\{
\begin{array}{l}
(i+l)\times (j_1+1)=n-m+1\\
l<i
\end{array} \right.
\end{equation}
From the lemma 3, (8) and (9) have $\frac {p-1}{2}$ solutions, respectively.

Thus, (5) has $p-1$ solutions. $\Box$\\

\end{document}